\documentclass[11pt]{article}
\usepackage{graphicx}
\usepackage{amsmath}
\usepackage{amsfonts}
\usepackage{amssymb}
\usepackage{makeidx}%
\newtheorem{theorem}{Theorem}
\newtheorem{lemma}{Lemma}
\newtheorem{definition}{Definition}

\newcommand{\R}  {\mathbb{R}}
\newcommand{\E}  {\mathbb{E}}
\renewcommand{\P}  {\mathbb{P}}

\newcommand{\F}  {\mathcal{F}}
\newtheorem{proposition}[theorem]{Proposition}
\newtheorem{remark}[theorem]{Remark}
\newcommand\unnumberedfootnote[1]{ %
        \let\temp=\thefootnote %
        \renewcommand{\thefootnote}{}%
        \footnote{#1}%
        \let\thefootnote=\temp%
        \addtocounter{footnote}{-1}}
        
\def \epf {\hbox{}\nobreak\hfill
\vrule width 2mm height 2mm depth 0mm
\par \goodbreak \smallskip}
\begin{document}
\title{From Brownian motion with a local time drift  \\to Feller's branching diffusion with logistic growth}
\author{ETIENNE PARDOUX\\[1ex]
  \emph{\small Universit\'e de Provence}\\
  e-mail: pardoux@cmi.univ-mrs.fr
  \\\\ANTON WAKOLBINGER \\
  \emph{\small Goethe--Universit\"at Frankfurt}
  \\ e-mail: wakolbinger@math.uni-frankfurt.de}
  \date{}
\maketitle
\unnumberedfootnote{\emph{AMS 2000 subject classification.} {\tt 60J70}
 (Primary) {\tt 60J55, 60J80, 60H10} (Secondary).}
\unnumberedfootnote{\emph{Keywords and phrases.}  {\tt Ray-Knight representation, local time, Feller branching with logistic growth, Brownian motion, local time drift, Girsanov transform.}}

\vspace{-1.5cm}
\begin{abstract} 
We give a new proof for a Ray-Knight representation of Feller's branching diffusion with logistic growth in terms of the local times of a reflected Brownian motion $H$  with a drift that is affine linear in the local time accumulated by $H$ at its current level. 
In \cite{LPW}, such a representation was obtained by an approximation through Harris paths that code the genealogies of particle systems. The present proof is purely in terms of stochastic analysis, and is inspired by previous work of Norris, Rogers and Williams \cite{NRW}.
\end{abstract}

\section{Introduction} The second one of the two classical Ray-Knight theorems (see e.g.\  \cite{RY} or  \cite{RW}) establishes a representation of Feller's branching diffusion in terms of reflected Brownian motion. In words, it may be stated as follows: Take a standard Brownian motion on $\mathbb R_+$ reflected at $0$, and stopped when the local time accumulated at $0$ reaches a value $x$. Then the (total) local time accumulated by the resulting path at ``height'' $t$, viewed as a process indexed by $t$, is a Feller branching diffusion obeying the SDE $dZ_t^x = 2\sqrt{Z_t^x} dW_t^x$ with $Z_0^x = x$. One way to interpret this  is to view the reflected Brownian path as an {\em exploration path} which codes the genealogy of a continuous state branching process (see e.g.  \cite{LG}): the local time of the exploration path at height $t$ measures the ``width of the genealogical forest'' at this level, or equivalently, the mass (or size) of the population that is alive at the time corresponding to this height. This mass is $Z_t$, the state of the branching process at time $t$. The exploration path is a concatenation of Brownian excursions, with each excursion corresponding to a continuum random tree in the sense of Aldous \cite {Al}. The independence of the excursions leads to independent increments of $(Z^x)$ when viewed as a process indexed by $x$. This mutual independence of offspring coming from different ancestors is also referred to as the {\em branching property}. 
%
%
%

Consider now, instead of a Feller branching diffusion, the same SDE with a logistic drift, namely
\begin{equation}\label{FelloG}
dZ^x_t=\left[\theta Z^x_t-\gamma\left(Z^x_t\right)^2\right]dt+2\sqrt{Z^x_t}dW_t,\quad Z^x_0=x,
\end{equation}
where $\theta, \gamma>0$ are given parameters which will be fixed throughout the paper, and $x>0$ is the initial population. 
Just like Feller's branching diffusion, also the solution of \eqref  {FelloG}, called {\em Feller's branching diffusion with logistic growth}, arises as the diffusion limit of discrete population models, but now with an interaction between individuals.  This interaction  can be thought of as a competition among individuals for resources, resulting in ``pairwise lethal fights'' (with intensity $\gamma$) that counteract the supercritical growth of the population. This population model and its diffusion limit \eqref{FelloG} have been studied  in \cite{AL}.

In this note we will specify an SDE for a process $(H_s)$, from which Feller's branching diffusion with logistic growth can be read off in the same way as Feller's branching diffusion is read off from reflected Brownian motion. Our proof will rely purely on stochastic analysis. Still, we give here some brief explanations on the underlying population model (for more illustrations and background we refer to our survey paper~\cite{PW}). 

The process $(H_s)$ will be reflected Brownian motion with a drift that depends on the local time $\ell$ accumulated at the current level $H_s$ up to time $s$. More specifically, the drift coefficient will be of the form $\theta/2 -\gamma \ell$.  One way to understand the form of the drift  is to see $(H_s)$ again as the exploration process of a forest of random real trees, and to think of an approximation in terms of piecewise linear, continuous processes with constant absolute slope (so-called {\em Harris paths}): the rate of minima (giving rise to new branches) is increased proportional to $\theta/2$, and the rate of maxima (describing deaths of branches) is increased proportional  to the number of individuals visited by the exploration process so far on the current level.
In our recent work \cite{LPW} we obtained the  process $H$ as the limit of exploration processes of discrete population models, and in this way provided a Ray-Knight representation of  \eqref{FelloG}. The discrete approach of \cite{LPW} gives a worthy insight into the result, since  the way how the genealogy is built and how the exploration process codes the genealogical tree of the population is readily understandable at the discrete level. 

The derivation presented in this note  does not rely on a discrete approximation, but directly exploits methods from stochastic analysis. We use ideas from previous work of Norris, Rogers and Williams
\cite{NRW} to which our attention was drawn after the completion of  \cite{LPW}  thanks to a hint of J-F Le~Gall.  In \cite{NRW} a generalization of the first Ray--Knight theorem for ``Brownian motions with a local time drift'' was provided for cases that include the drift appearing in the SDE~\eqref{SDEexplor}.

We also extend the Ray-Knight representation of \eqref{FelloG} by establishing an equality between laws of random fields (random functions of time $t$ and ancestral mass $x$). In Section \ref{coupling} we introduce Feller's branching diffusion with logistic growth as a random field $\{Z^x_t,\ t,x\ge0\}$. This is a natural set-up for the formulation of our main result, which is given in Section \ref{main} and whose proof is contained in Section \ref{proof}. The last section gives two remarks concerning a possible shortcut in the proof of the Theorem, and a general version of the second Ray-Knight theorem in the framework of  \cite{NRW}.
\section{A coupling  over the ancestral masses}\label{coupling}
In this section we define a random field $\{Z^x_t,\ t, x\ge0\}$ such that for any $x\ge0$, $Z^x := \{Z^x_t,\ t\ge0\}$ is a Feller branching diffusion with logistic growth and ancestral mass $x$, for any $t \ge0$, $x \mapsto Z^x_t$ is non-decreasing and
$x \mapsto \{Z^x_t, \, t \ge 0\}$ is a (path-valued) Markov process which will be specified below..

To this purpose we define a family of transition probabilities $\mathbf P_x, x \ge 0$, on $E$, where $E=C^{c}(\mathbb R_+,\mathbb R_+)$  is the set of continuous mappings from $\mathbb R_+$ to  $\mathbb R_+$ with compact support. Here and below, $\mathbb{R}_+=[0,+\infty)$. For $x > 0$ and $z \in C^{c}(\mathbb R_+,\mathbb R_+)$, let $\mathbf P_x(z,\cdot)$ be the distribution of $z + Z^{z,x}$, where $Z^{z,x}$ solves
\begin{equation}\label{fellogz}
Z^{x,z}_t=x + \int_0^t Z^{x,z}_u (\theta-\gamma [ Z^{x,z}_u +2z(u)])du+2\int_0^t\sqrt{Z^{x,z}_u}dW_u,
\end{equation}
with $W$ being a standard Brownian motion. The equality  $\mathbf P_x(z,E)=1$ is valid because $Z^{x,0}$ (and a fortiori $Z^{x,z}$) a.s. hits zero in finite time (for a proof of this fact see e.g. \cite{AL}).
Before we show in Lemma \ref{ChKo} that the transition probabilities $\mathbf P_x, x \ge 0$, indeed specify a random field $\{Z_t^x\}$, let us briefly motivate the form of \eqref {fellogz} in the light of the interpretation given in the introduction.

Consider the evolution of the progeny of a sum of ancestral masses $z(0)+x$. The offspring of $z(0)$ is assumed to follow (in an autonomous way) the dynamics of a Feller branching with logistic drift,  the path $z(t)$ stands for a realization of this offspring. The progeny $(Z_t^{x,z})$ of the additional ancestral mass $x$ does not evolve independently of the offspring of $x$, but  experiences an additional pressure coming from the given $z(t)$, resulting in the negative drift $-2 \gamma z(t)Z_t^{x,z} $. In a finite population approximation with $k+\ell$ ancestors, this means  
that the descendants of the $\ell$ ``additional'' ancestors suffer from the competition with those of the ``first'' $k$, while the descendants of the ``first'' $k$ do not feel the presence of the descendants of the $\ell$ ``additional'' ancestors. Recall from the introduction that in an individual-based model the competition leading to the negative quadratic drift in the populaiton size  is modelled by pairwise fights between the individuals. If we think of the individuals being arranged in a linear order ``from left to right'', where this order is passed on to the individual's offspring, then in our convention the pairwise fights  are always be won by the individual to the left, resulting in the death of the individual to the right.  
\begin{lemma}\label{ChKo}
The family  $\mathbf P_x, x \ge 0$, satisfies the Chapman-Kolmogorov relations.
\end{lemma}
 {\sc Proof:}
Observe that conditioned on $Z^{x,z}$, the random path $V:= Z^{y, z+Z^{x,z}}$ solves 
\begin{equation}\label{fellogV}
V_t=y + \int_0^tV_u (\theta-\gamma [V_u +2(z(u)+Z^{x,z}_u)])du+2\int_0^t\sqrt{V_u}dW'_u
\end{equation}
with $W'$ being a standard Brownian motion (independent of $W$). Note that in the case $\gamma=0$, the two processes $Z^{x,z}$ and $V$ are independent, as they should.  Now $Z^{x,z} + V$ satisfies
\begin{align*}
Z^{x,z}_t + V_t = x+y+\int_0^t(Z^{x,z}_u + V_u)(\theta-\gamma[Z^{x,z}_u + V_u+2z(u)])du\\ +2\int_0^t\sqrt{Z^{x,z}_u}dW_u+
2\int_0^t\sqrt{V_u}dW'_u\, .
\end{align*}
This shows that $z+Z^{x,z}+ V $ has distribution $\mathbf P_{x+y}(z,\cdot)$, as required.
 Indeed, since the two Brownian motions $W$ and $W'$ are independent, 
$$\langle2\int_0^\cdot \sqrt{Z^{x,z}_u}dW_u+2\int_0^\cdot \sqrt{V_u}dW'_u\rangle_t=
4\int_0^t(Z^{x,z}_u+V_u)du,$$ and consequently, from a well--known martingale representation theorem, there exists a third standard Brownain motion $W''_t$ such that
$$2\int_0^t \sqrt{Z^{x,z}_u}dW_u+2\int_0^t \sqrt{V_u}dW'_u=2\int_0^t\sqrt{Z^{x,z}_u+V_u}dW''_u. \quad \mbox{\epf}$$
\begin{definition}
Let  $\{Z^x\}_{x\ge 0}$ be the  $C^c(\mathbb R_+, \mathbb R_+)$-valued Markov process with transition semigroup $(\mathbf P_x)$, starting from $Z^0 \equiv 0$, the null trajectory.
\end{definition}
\begin{remark} For each $x > 0$, $Z^x$ solves the SDE 
\begin{equation}\label{fellog}
dZ^x_t=\left[\theta Z^x_t-\gamma(Z^x_t)^2\right]dt+2\sqrt{Z^x_t}dW_t^x,\ Z^x_0=x,
\end{equation}
where $\{W_t^x,\ t\ge0\}$ is a standard Brownian motion. Since for $x, y > 0$ the increment  $Z^{x+y}-Z^x$ is driven by a Brownian motion independent of that driving $Z^x$, we have 
$d\langle Z^x, Z^{x+y}\rangle_t = d\langle Z^x, Z^x\rangle_t = Z_t^x\,dt $
and conseqently
$d\langle W^x,W^{x+y}\rangle_t=\sqrt{{Z^x_t}/{Z^{x+y}_t}}\,dt,\quad \text{with the convention }\frac{0}{0}=0.$
\end{remark}

\section{A Ray-Knight representation}\label{main}
Consider the following SDE driven by standard Brownian motion~$B$
\begin{equation}\label{SDEexplor}
H_s=B_s+\frac{1}{2}L_s(0)+\frac{\theta}{2}s-\gamma\int_0^sL_r(H_r)dr,\ s\ge0,
\end{equation}
Here and everywhere below, $\{L_s(t),\ s\ge0,\ t\ge0\}$ denotes the local time of the process $\{H_s,\ s\ge0\}$
accumulated up to time $s$ at level $t$. Proposition \ref{GirExp}, stated and proved in the next section,
will ensure (by specializing it to the case $z\equiv 0$) that  equation \eqref{SDEexplor} has a unique
 weak solution, which we assume to be defined on some probability space $(\Omega,\F,\P)$.  
 
 Define for any $x>0$ the stopping time
$$S_x=\inf\{s>0,\  L_s(0)>x\},$$
and let $\{Z^x_t,\ x,t\ge 0\}$ denote the random field constructed in Section \ref{coupling}.\\\\
Our main result is the\\\\
{\bf Theorem} {\em
The two random fields $\{L_{S_x}(t),\ t, x\ge0\}$ and $\{Z^x_t,\ t, x\ge0\}$
have the same law.}

\section{Proof of the Theorem }\label{proof}
To prepare for the proof of the Theorem, we first fix a $z \in C^{c}(\mathbb R_+,\mathbb R_+)$ and consider the SDE
\begin{equation}\label{SDEexplor-z}
H^z_s=B_s+\frac{1}{2}L^z_s(0)+\frac{\theta}{2}s-\gamma\int_0^s\{z(H^z_r)+ L^z_r(H^z_r)\}dr,\ s\ge0,
\end{equation}
where $L^z$ stands for the local time of $H^z$. We will prove in Subsection~\ref{ss4.1} 
\begin{proposition}\label{GirExp}
The SDE \eqref{SDEexplor-z} has a unique weak solution.
\end{proposition}  

Suppressing the superscript $z$, define for any $x>0$ the stopping time
\begin{align} \label{Sz} S_x=\inf\{s>0,\  L^z_s(0)>x\}.\end{align}
 The main step in the proof of the Theorem will be to show
 \begin{proposition}\label{RNz} For  $x>0$ and $z\in C^c(\R_+,\R_+)$ let $\{Z^{x,z}_t,\ t\ge0\}$ be the solution of \eqref{fellogz}.
Then the two processes $\{L^z_{S_x}(t),\ t\ge0\}$ and $\{Z^{x,z}_t,\ t\ge0\}$
 have the same law.
 \end{proposition}

\subsection{Proof of Proposition \ref{GirExp}}\label{ss4.1}
Let $H$ denote Brownian motion reflected above 0, i. e.
$$H_s=B_s+\frac{1}{2}L_s(0),$$
where  $B$ is a $\F_s$--standard Brownian motion defined on a probability space $(\Omega,\F,\P)$, 
with $\F=\F_\infty$, and $L$ is the semimartingale local time of $H$.  
Let $$G_s=\exp\left(M_s-\frac{1}{2}\langle M\rangle_s\right),\quad s\ge0,$$
with $M_s := \int_0^s\left\{\frac{\theta}{2}-\gamma\left[z(H_r)+L_r(H_r)\right]\right\}dB_r$. (Recall that $z \in C^{c}(\mathbb R_+,\mathbb R_+)$ is fixed.) The condition 
\begin{equation}\label{Girs}
 \mathbb{E}( G_s) =1, \quad \forall s\ge 0
 \end{equation}
 is sufficient  for the local martingale $\{G_s,\ s\ge0\}$ to be  a martingale. In that case,
 there exists a new probability measure $\tilde{\P}$ on $(\Omega,\F)$ such that for all $s\ge0$, 
 $$\frac{d\tilde{\P}|_{\F_s}}{d\P|_{\F_s}}=G_s,$$
 and it follows from Girsanov's theorem (see e. g. Theorem VIII 1.4 in \cite{RY}) that
 \begin{equation}\label{Btilde}
\tilde {B}_s:=B_s-\int_0^s\left\{\frac{\theta}{2}-\gamma\left[z(H_r)+L_r(H_r)\right]\right\}dr, \qquad s\ge 0,
\end{equation}
 is a standard Brownian motion.  (Note that this does not require that $\tilde{\P}$ be absolute continuity with respect to $\P$ on $\F$.) Hence
existence of a weak solution to \eqref{SDEexplor-z} follows from \eqref{Girs}, which in turn
 (see Theorem 1.1, chapter 7, page 152, in \cite{Fr}) follows if we can ensure that for each $s>0$ there exists constants $a>0$ such that
 \begin{equation}\label{suffgir}
  \sup_{0\le r\le s}\mathbb E \exp(aR_r) <\infty,
  \end{equation}
where $R_r=\left|\frac{\theta}{2}-\gamma\left[z(H_r)+L_r(H_r)\right]\right|^2$. Since $z$ is bounded, the inequality \eqref{suffgir} is immediate from the following
 \begin{lemma}
Let $H$ be a Brownian motion   on $\mathbb R_+$ reflected at the origin.
Then for all $s>0$ there exists $\alpha = \alpha(s) >0$ such that   
$$\sup_{0\le r\le s}\mathbb E\left(\exp (\alpha L_r(H_r)^2)\right)<\infty.
$$
\end{lemma} 
\noindent
 {\sc Proof:} Together with a simple scaling argument and a desintegration with respect to $H_r$, this is immediate from the following
 \begin{lemma}\label{monloctime}
Let $\beta$ be a standard Brownian motion starting at $0$, and denote by $L_1(y)$ the local time accumulated by $|\beta|$ at position $y$ up to time $1$. There exist  constants $a >0$ and $c >0$ (not depending on $y$) such that for almost all $y \ge 0$
 \begin{equation} \label{expmom}\mathbb E[e^{aL_1(y)^2}|\, |\beta_1| = y] \le c.
  \end{equation}
\end{lemma} \noindent
 {\sc Proof:} By symmetry the l.h.s. of \eqref {expmom} a.s. equals $\mathbb E[e^{aL_1(y)^2}|\, \beta_1 = y].$   Writing $\mathbb P^y$ for the probability measure of a Brownian bridge from the origin at time $0$ to position $y$ at time $1$, and $K_1(a)$ for the local time accumulated up to time $1$ at position $a$, we thus have to show the inequality 
\begin{equation} \label{BB} E[e^{a(K_1(y) +  K_1(-y))^2}] \le c
\end{equation}
for suitable constants $a$ and $c$.
By  the Cauchy-Schwarz inequality, the l.h.s. of \eqref{BB} is bounded by
\begin{equation}\label{CS}
\left(\mathbb E^y\left[e^{4a K_1(y)^2}\right]\right)^{1/2} \left(\mathbb E^y\left[e^{4a K_1(-y)^2}\right]\right)^{1/2} .
\end{equation}
To estimate the first factor, we desintegrate with respect to the time $U$ at which the path of the Brownian bridge first hits the  level $y$. Conditioned under $\{U=u\}$, the part before time $u$ does not contribute to the local time at $y$, and the second part is (by the stong Markov property) a Brownian bridge from $y$ to $y$ over a time interval of length $1-u$, hence (by scaling) the distribution of its local time at $y$ equals the distribution of  $\sqrt{1-u}\, K_1(0)$ under $\mathbb P^0$. We therefore obtain
\begin{equation}\label{desint}\mathbb E^y\left[e^{4a K_1(y)^2}| U=u\right] = \mathbb E^0\left[e^{4a (1-u) K_1(0)^2}\right]\le \mathbb E^0\left[e^{4a K_1(0)^2}\right]
\end{equation}
To estimate the second factor, we desintegrate with respect to the times $(U_1, U_2)$ at which the path of the Brownian bridge hits the position $-y$ for the first resp. the last time (on the event that it hits this position at all). Again by the strong Markov property and scaling, the distribution $K_1(-y)$ under $\mathbb P^y[.| \{U_1=u_1, U_2=u_2\}]$  equals the distribution of  $\sqrt{u_2-u_1}\, K_1(0)$ under $\mathbb P^0$. This leads to an estimate analogous to \eqref{desint}, and allows to conclude  
\begin{equation}\label{compa}
\mathbb E^y\left[e^{4a K_1(y)^2}\right]\le \mathbb E^0\left[e^{4a K_1(0)^2}\right] , \quad \mathbb E^y\left[e^{4a K_1(-y)^2}\right]\le \mathbb E^0\left[e^{4a K_1(0)^2}\right].
\end{equation}
By a result due to L\'evy (see formula (11) in \cite{Pi}),   $K_1(0)$ has under $\mathbb P^0$ a  Raleigh distribution, i.e.
$$\mathbb P^0(K_1(0)   > \ell ) = e^{-\frac 12 \ell^2}.$$
This means that $K_1^2(0)$ is  exponentially distributed, and hence, for suitably small $\delta>0$,  $\mathbb E^0\left[e^{\delta K_1(0)^2}\right]$ is finite. Now \eqref{expmom} follows from \eqref{CS} and \eqref{compa}. 
\epf

 So far we have proved existence of a weak solution to \eqref{SDEexplor-z}. Weak uniqueness is easier to prove, since uniqueness is a local property. Let $H$ be a solution to equation \eqref{SDEexplor-z}, and for all $n\ge1$ let $T_n$ denote the stopping time
 $$T_n := \inf \{r > 0 : L_r(H_r) > n\}.$$
 By a Girsanov transformation we can change the measure $\mathbb P$ into a measure $\bar {\mathbb P}$ under which, for all $n \in \mathbb N$,  the restriction of the process $H$ to the interval $[0,n\wedge T_n]$ is standard Brownian motion reflected above 0.
 Since $\mathbb P$ and $\bar {\mathbb P}$ are mutually absolutely continuous, the law of $\{H_{s\wedge n\wedge  T_n},\ s\ge0\}$ under $\mathbb P$ is uniquely determined, for each $n\ge1$. Uniqueness of the law of $H$ solution of \eqref{SDEexplor-z} then
 follows, since $T_n\to\infty$ a. s. as $n\to\infty$.

\subsection{Proof of Proposition \ref{RNz}}
As a by-product of our proof, we will see that the stopping time $S_x$ defined in \eqref{Sz} has finite expectation. A more direct argument for this would make the proof of Proposition \ref{RNz} even shorter, see the discussion in Subsection \ref{Disc}. Since we have not been able to prove this directly, we circumvent this by reflecting the process $H^z$ below the level $K$, and then let $K$ tend to $\infty$. 

To be specific,
for $K>0$, let $H^K$ be the solution of the SDE
\begin{align} \label{refl} H^K_s= B_s+\frac{1}{2}L^K_s(0)-\frac{1}{2}L^K_s(K^-),\ s\ge0,
\end{align}
where $L^K$ denotes the local time of $H^K$ and $B$ is standard Brownian motion defined on the probability space  $(\Omega,\F,\P)$. In other words, $H^K$ is Brownian motion reflected inside the interval $[0,K]$.

Let us first note that if we define
\begin{equation}\label{SKdef}
S^K_x=\inf\{s>0,\  L^K_s(0)>x\},
\end{equation}
the next result follows readily from Lemma 2.1 in Delmas \cite{JFD}:
\begin{lemma}\label{leJFD}
For any $K>0$ the processes $\{L_{S_x}Ñ(t),\ 0\le t\le K\}$ and $\{L^K_{S^K_x}Ñ(t),$ $\ 0\le t\le K\}$
have the same distribution.
\end{lemma}
(The intutive explanation of this lemma is as follows: Consider an arbitrary level $K>0$. The law of Brownian motion reflected in $[0,K]$ equals the law of Brownian motion reflected above 0, from which we chop out the pieces of trajectories which exceed the level $K$.)

We next define the martingale
$$M^K_s=\int_0^s \left[\frac{\theta}{2}- \gamma\{z(H^K_r)+ L_r^K(H^K_r)\}\right]dB_r.$$
The same arguments as those in the proof of Proposition \ref{GirExp} show here also that for all $s>0$,
$$\E\exp\left( M^K_s-\frac{1}{2}\langle M^K\rangle_s\right)=1.$$
Therefore there exists a probability measure $\tilde{\P}^K$ such that for all $s>0$, 
$$\frac{d\tilde{\P}^K}{d\P}\left|_{\F_s}=\exp\left( M^K_s-\frac{1}{2}\langle M^K\rangle_s\right).\right.$$
From Girsanov's theorem, under $\tilde{\P}^K$, $H^K$ is a solution of the reflected SDE
\begin{equation}\label{explorz}
H^K_s= B_s+\frac{\theta}{2}s
-\gamma\int_0^s[z(H^K_r)+L^K_r(H^K_r)]dr+\frac{1}{2}L^K_s(0)-\frac{1}{2}L^K_s(K^-),\ s\ge0.
\end{equation}
We will require
\begin{lemma}\label{SKxlem}
 \begin{equation*}\label{SKx}
\tilde{\E}^K[S^K_x]< \infty.
\end{equation*}
\end{lemma}
{\sc Proof:}
We will prove this by a comparison argument. To this end let, under $\tilde{\P}^K$,  $\bar H^K$ be the solution of 
 \begin{equation}\label{explorsuperc}
\bar H^K_s= B_s+\frac{\theta}{2}s
+\frac{1}{2}\bar L^K_s(0)-\frac{1}{2}\bar L^K_s(K^-),\ s\ge0,
\end{equation}
with $\bar L^K$ denoting the local time of $\bar H$; in other words, $\bar H^K$ is a Brownion motion with constant upward drift $\theta/2$, reflected above $0$ and below $K$. Obviously, for all $n=1,2, \ldots$ , $a\in [0,K]$ and $\varepsilon > 0$,
$$\tilde \P (\bar L^K_{n+1}(0)-\bar L^K_n(0) \ge \varepsilon | \bar H^K_n = a) \ge \tilde \P (\bar  L^K_{n+1}(0)-\bar L^K_n(0)  \ge \varepsilon | \bar H^K_n = K) := \delta,$$
where $\delta$ depends on $\varepsilon$ and $K$ but not on $n$ and $a$, and is positive at least for sufficiently small $\varepsilon$. Choosing this $\varepsilon$ and $\delta$, we see that the probability that $\bar L^K(0)$ increases in the time interval $[n,n+1]$ by at least $\varepsilon$, bounded from below by $\delta$, independent of the past. This implies 
 \begin{equation}\label{SKbar}
\tilde{\E}^K[\bar S^K_x]< \infty,
\end{equation}
where $\bar S^K_x$ is defined by \eqref{SKdef}, there with $L_s^K(0)$ replaced by $\bar L_s^K(0)$.
From a classical 
comparison theorem for SDEs, see e.g. Theorem 3.7, chapter IX of \cite{RY}, we conclude that 
$H^K_s\le \bar{H}^K_s$, $s\ge 0$, a.s. 
 This implies that 
$$\bar{L}^K_{s}(0)\le L^K_{s}(0), \ s\ge 0,  \mbox { a.s. }$$
Consequently, $\bar S^K_x \ge S^K_x$ a.s.; hence the assertion follows from  \eqref{SKbar}.
\epf
The next subsection will be devoted to the proof of 
\begin{proposition}\label{RNK}
For any $K>0$, the process $\{L^K_{S^K_x}(t),\ t\ge0\}$ is under $\tilde{\P}^K$ a solution of equation \eqref{fellogz},
killed at time $K$.
\end{proposition}
From this together with Lemma~\ref{leJFD}, 
Proposition \ref{RNz}  is immediate.

\subsection{Proof of Proposition \ref{RNK}}
In this section, $x>0$ and $K>0$ are fixed. We work under  $\tilde{\P}^K$  and take advantage of some of the techniques  
from \cite{NRW}.
\\\\
 Tanaka's formula yields for any $r \ge 0$ and $0\le t<K$ the identity 
 $$(H^K_r-t)^- = (-t)^- - \int_0^{r}{\bf1}_{\{H^K_s\le t\}}dH^K_s+\frac{1}{2}
 L^K_{r}(t).$$
 With $r:=S_x^K$ (which is finite  $\tilde{\P}^K$-a.s. due to Lemma \ref{SKxlem}) this yields
\begin{equation} \label{loctime} L^K_{S^K_x}(t) = 2 \int_0^{S^K_x}{\bf1}_{\{H^K_s\le t\}}dH^K_s\ .
\end{equation}
 Plugging \eqref{explorz} into \eqref{loctime} we arrive at
 \begin{eqnarray}\nonumber
 L^K_{S^K_x}(t) &=& x+ 2 \int_0^{S^K_x}{\bf1}_{\{H^K_s\le t\}}dB_s
 \\ \label{loctimedyn} &+& \int_0^{S^K_x}{\bf1}_{\{H^K_s\le t\}}\left( \theta -2\gamma\{z(H^K_s)+ L^K_s(H^K_s)\}\right)ds.
 \end{eqnarray}
 It follows from the occupation times formula (see e.g. \cite{RY} VI.1.7) that
\begin{equation}\label{occ1}
\int_0^{S^K_x}{\bf1}_{\{H^K_s\le t\}}\left( \theta -2\gamma \, z(H^K_s)\right)ds
 =\int_0^t \left( \theta -2\gamma \, z(u)\right)L^K_{S^K_x}(u)du,
 \end{equation}
 while from a generalization of the same formula (see Exercise 1.15 in Chapter VI of \cite{RY}) we have
 \begin{align} \nonumber
 2\gamma\int_0^{S^K_x}{\bf1}_{\{H^K_s\le t\}}L^K_s(H^K_s)ds&=
 2\gamma\int_0^t\int_0^{S^K_x}L^K_s(u)dL^K_s(u)du\\ \label{occ2}
 &=\gamma\int_0^t\left(L^K_{S^K_x}(u)\right)^2du.
 \end{align}
Let us now  abbreviate 
\begin{equation}\label{Nt1} N_t :=2 \int_0^{S^K_x}{\bf1}_{\{H^K_s\le t\}}dB_s, \quad 0\le t \le K.
\end{equation}
In order to check that this is a martingale with the appropriate quadratic variation, we define, following \cite{NRW},  for all  $0\le t\le K$ and $s \ge 0$
 \begin{align*}
A(s,t):=\int_0^s{\bf1}_{\{H^K_r\le t\}}dr&,\quad \tau(r,t):=\inf\{s:\ A(s,t)>r\},\\
J(s,t):=\int_0^s{\bf1}_{\{H^K_r\le t\}}dB_r,&\quad \xi(r,t):=J(\tau(r,t),t).
\end{align*}
For fixed $t$, the process $\xi(.,t)$ as a Brownian motion, arising through a time change from the continuous martingale  $J(.,t)$.
Write $\mathcal F(.,t)$ for the  filtration generated by $\xi(.,t)$, and  $\mathcal E_t := \mathcal F(\infty,t)$ for the $\sigma$-algebra generated by the  $\mathcal F(s,t)$, $0\le s < \infty$. With these slight modifications of the definitions given in \cite{NRW} p. 273, we can carry over all the steps in the proof of \cite{NRW}, Theorem 1, to our situation. We will explain here the main ideas and a few details. 

A crucial observation  is that every bounded $\mathcal E_t$--measurable random variable $F$ can be represented as an It\^o integral
\begin{equation}\label{testF}
F= \mathbb E [F] + \int_0^\infty v_r d\xi(r,t) = \mathbb E[F] + \int_0^\infty v_{A(s,t)} \mathbf 1_{\{H_s \le t\}} dB_s
\end{equation}
for some $\mathcal F(.,t)$--predictable process $v$ such that $\mathbb E\int_0^\infty v_r^2 dr < \infty$.

Let $u_s := 2\cdot \mathbf 1_{\{0\le s \le S_x^K\}}$. This process is predictable, and Lemma~\ref{SKxlem}
implies that $\mathbb{E}\int_0^\infty u_s^2ds<\infty$. Moreover for each $t>0$, the process
$$\tilde{u}(s,t),\ s\ge0\ \text{is }\mathcal F(.,t)\text{--predictable},$$
since $\tilde{u}(s,t)=2\cdot \mathbf 1_{\{0\le s \le A(S_x^K,t)\}}$, and $A(S_x^K,t)$ is a 
$\mathcal F(.,t)$--stopping time.

Consequently (recall \eqref{Nt1})
$$N_t=\int_0^\infty u_s{\bf1}_{\{H^K_s\le t\}}dB_s=\int_0^\infty \tilde{u}(s,t)d\xi(s,t)$$
is $\mathcal E_t$--measurable, as well as
$$C_t=\int_0^\infty u_s^2{\bf1}_{\{H^K_s\le t\}}ds=\int_0^\infty \tilde{u}(s,t)^2ds.$$

Now writing $N_{t+h}-N_t$ as the integral $$\int_0^\infty u_s\mathbf 1_{\{t<H_s\le t+h\}}dB_s,$$ we see from \eqref{testF} that 
$\mathbb E[(N_{t+h}-N_t)F] = 0$, in other words, $\{N_t,\ 0\le t\le K\}$ is an $(\mathcal{E}_t)$--martingale. 

Applying It\^o's formula to the process $s \mapsto \left(\int_0^s u_r\mathbf 1_{\{t<H_r\le t+h\}}dB_r\right)^2$, \, $s \ge 0$, we obtain as in \cite{NRW} that 
$$\mathbb E[(N_{t+h}-N_t)^2F] = \mathbb E[\int_0^\infty u^2_s \mathbf 1_{\{t<H_s\le t+h\}}ds \cdot F ]=  4\mathbb E[\int_0^{S^K_x}\mathbf 1_{\{t<H_s\le t+h\}}ds \cdot F ],$$
which reveals 
 the quadratic variation of $(N_t)$ as
$
 \langle N\rangle_t=4\int_0^{S^K_x}{\bf1}_{\{H^K_s\le t\}}\, ds $.
Again by the occupation times formula, this equals $2\int_0^t L^K_{S^K_x}(u)\, du$. Consequently, there exists a Brownian motion $\{W_t,\ t\ge0\}$ such that
\begin{equation} \label{Nt} N_t = 2\int_0^t  \sqrt{ L^K_{S^K_x}(u)} dW_u, \quad 0\le t \le K.
\end{equation}
 The proof of Proposition \ref{RNz} is now completed by combining \eqref{loctimedyn}, \eqref{occ1}, \eqref{occ2} \eqref{Nt1} and \eqref{Nt}.
 \epf

\subsection{Completion of the proof of the Theorem}
It follows from the description  of the law of $\{Z^x_t,\ t\ge0\}_{x\ge0}$ made in Section~\ref{coupling} that
$Z$ is Markov (as a process indexed by $x$, with values in the set  of continuous paths from
$\R_+$ into $\R_+$ with compact support). The fact that
 $\{L_{S_x}(t),\ t\ge0\}_{\{x\ge0\}}$ enjoys the same property follows from the fact that the process
 $H^x_r:=H_{S_x+r}$ solves the SDE \eqref{SDEexplor-z} with $z(t)=L_{S_x}(t)$ and a Brownian motion $B$ which,
 from the strong Markov property of Brownian motion, is independent of  $\{L_{S_x}(t),\ t\ge0\}$.
 
  Hence it suffices to prove that for any $0\le x<x+y$, the conditional law of 
$L_{S_{x+y}}(\cdot)$ given $L_{S_x}(\cdot)$ equals that of $Z^{x+y}_\cdot$, given $Z^x_\cdot$. Conditioned upon
$L_{S_x}(\cdot)=z(\cdot)$, $L_{S_{x+y}}(\cdot)-L_{S_x}(\cdot)$ is the collection of local times accumulated by the solution of 
\eqref{SDEexplor-z} up to time $S_{y}$, i. e. it has the law of the process $\{L^z_{S_{y}}(t),\ t\ge0\}$, while conditionally upon  $Z^x_\cdot=z(\cdot)$, the law of $Z^{x+y}_\cdot-Z^x_\cdot$
is that of $Z^{y,z}$, solution of equation \eqref{fellogz}. Thus, the assertion of the Theorem follows from Proposition \ref{RNz}.

\section{Concluding remarks}
\subsection{A possible shortcut in the proof of Proposition \ref{RNz} }\label{Disc}
As a direct consequence of Proposition \ref{RNz} and the occupation times formula, 
the stopping time $S^x$ defined by \eqref{Sz} obeys
$$S_x\stackrel d = \int_0^\infty Z^{x,z}_t dt,$$
where $Z^{x,z}$ is the solution of \eqref{fellogz}.
This together with a representation of  $\int_0^\infty Z^{x,0}_t dt$ as the random time at which an Ornstein-Uhlenbeck process first hits $0$  (see \cite{AL}) proves
\begin{lemma}\label{Sxfinite}
For any $x>0$, the stopping time $S_x$ defined in \eqref{Sz} has finite expectation.
\end{lemma}
If we could prove Lemma \ref{Sxfinite} directly from the SDE \eqref{SDEexplor-z}, then we could simplify our proof of Proposition \ref{RNz}, avoiding the 
reflection below the arbitrary level $K$. Here is an attempt of a direct intuitive explanation why  Lemma \ref{Sxfinite} holds.
 While climbing up, the Brownian motion with positive drift $\theta/2$ accumulates local time at various levels. Sooner or later, it accumulates so much local time around some level in $\mathbb{R}_+$ that the process $H$ governed by \eqref{SDEexplor} starts to go down. It then continues to accumulate local time at various levels, and goes back to zero. After reflection at zero, the next excursions will have already a stronger drift downwards that awaits $H$. Remarkably, the recurrence  of $H$ to the state $0$ holds independently of the relative constellations of the positive parameters $\theta$ and $\gamma$. 
 
 \subsection{A second Ray-Knight theorem for Brownian motion with a local time drift}
 The equation \eqref{SDEexplor-z} is of the form
 \begin{equation}\label{SDEexplor-gen}
H_s=B_s+\frac{1}{2}L_s(0)+\int_0^s g(H_r, L_r(H_r))dr,\ s\ge0,
\end{equation}
The proof of Proposition \ref{RNz} shows that $\{L_{S_x}(t), t \ge 0\}$ satisfies the SDE 
\begin{equation}\label{fellog-gen}
Z_t=x + \int_0^t f(u, Z_u) du + 2\int_0^t\sqrt{Z_u}dW_u
\end{equation}
with $f(t,\ell)= \int_0^\ell g(t,y)dy$, provided $g$ is such that \eqref{SDEexplor-gen} and \eqref{fellog-gen} have unique weak solutions which arise via Girsanov transformations from the  distributions with $g\equiv 0$, and provided $S_x=\inf\{s>0,\  L_s(0)>x\} $ is finite a.s. This more general problem will be the object of a forthcoming paper, where we will in particular  make precise the  interaction inside the population, at the discrete population level, leading to the continuous limit \eqref{fellog-gen}.
\\\\
{\bf Acknowledgement:} We thank Jean-Fran\c cois Le Gall for drawing our attention to the paper \cite{NRW}, Yueyun Hu for helping us streamline the proof of Proposition  \ref{RNK}, and a referee for a careful reading of a first version that led to an improved presentation.

\end{document}